\documentclass[oneside,12pt]{amsart}
\newtheorem{thm}{Theorem}
\newtheorem{lem}{Lemma}
\newtheorem{prop}{Proposition}
\input amssym.def
\input amssym.tex

\newcommand{\beq}{\begin{equation}{}}
\newcommand{\bge}{\begin{equation}{}}

\newcommand{\CCC}{{\Bbb C}}

\newcommand{\eed}{\ :=\ }
\newcommand{\eee}{\ =\ }

\newcommand{\fl}{{\par\noindent}}

\newcommand{\goin}{{\ \to \infty}}

\newcommand{\inv}{^{-1}}

\newcommand{\lee}{{\ \le \ }}

\newcommand{\MM}{{\mathcal M}}
\newcommand{\NN}{{\Bbb N}}
\newcommand{\pard}{\partial}

\newcommand{\pf}{\fl{{\bf Proof}. }}
\newcommand{\PP}{{\mathcal P}}

\newcommand{\RR}{{\Bbb R}}

\newcommand{\sumi}{\sum^{\infty}}

\newcommand{\ZZ}{{\Bbb Z}}

\newcommand{\bbe}{\vskip .15pc\hangindent=1pc\hangafter=1\noindent}
\def\ffdr{{\mathcal F}_{\delta,r}}
\def\ffdv{{\mathcal F}_{\delta,v}}
\def\ggdr{{\mathcal G}_{\delta,r}}

\def\ffdro{{\mathcal F}_{\delta,r}^{(1)}}
\def\ffdrt{{\mathcal F}_{\delta,r}^{(2)}}
\def\ffdrj{{\mathcal F}_{\delta,r}^{(j)}}
\def\ffdvj{{\mathcal F}_{\delta,v}^{(j)}}
\def\ggdrj{{\mathcal G}_{\delta,r}^{(j)}}
\def\nrmdrstj{{\|^{\ast,j}_{\delta,r}}}
\def\anormdstj{{$\|\cdot\nrmdrstj$}}
\def\mm{{\MM\MM}}
\def\mmr{{{\MM\MM}_r}}
\def\mmtr{{{\MM\MM}_{2r}}}

\def\wwdelt{{{\mathcal W}_{\delta}}}
\begin{document}
\centerline{Some classes of rational functions and related Banach spaces }
\vskip .1in
\centerline{by R.\ M.\ Dudley, Sergiy Sidenko, Zuoqin Wang, and Fangyun Yang}
\vskip .1in
\centerline{Department of Mathematics, Massachusetts Institute of Technology}
\vskip .1in
\centerline{\footnotesize All authors partially supported by National
Science Foundation Grant DMS-0504859}
\vskip .1in
Throughout this paper, $d$ and $r$ will be positive integers
and $0<\delta<1$. These numbers will be arbitrary with the
given properties unless they are further specified.

For a matrix or vector $X$ let $X'$ denote its transpose. Thus if
$x=(x_1,...,x_d)'\in\RR^d$ is a column vector and $A$ is a $d\times d$
matrix, then $x'Ax$ gives the quadratic form defined by the matrix
$A$, $\sum_{i,j=1}^n A_{ij}x_ix_j$. 
Let $\PP_d$ denote the set of
symmetric $d\times d$ positive definite matrices. 
Let $\|A\|$ be the
usual operator norm of a matrix $A$, $\|A\|\eed\sup\{|Ax|:\ |x|=1\}$,
where $|\cdot|$ is the usual Euclidean norm on $\RR^d$.
Recall that $\NN$ is the set of nonnegative integers. 
Let $\mmr\eed {\MM\MM}_{r,d}$ be the set of monic monomials from 
$\RR^d$ into $\RR$ of degree 
$r$, namely the set of all functions 
$g(x)=\Pi^d_{i=1}x_i^{n_i}$ with
$n_i\in\NN$ and $\sum^d_{i=1}n_i
=r$. 
Let
$\wwdelt\eed{\mathcal W}_{\delta,d}\eed\{C\in\PP_d:\ 
\|C\|<1/\delta,\ \|C\inv\|<1/\delta\}$.
Let
$$
\ffdr\eed{\mathcal F}_{\delta,r,d}\eed
\Big\{f:\ \RR^d\to\RR,\ f(x)\equiv  
g(x)/\Pi_{s=1}^{r}
(1 + x'C_s x), 
$$
$$
\textrm{where \ } 
g\in\mm_{2r},\textrm{ and for \ }s=1,...,r,\ 
C_s\in\wwdelt \Big\}.
$$
Let $\ggdr:=\bigcup_{v=1}^r\ffdv$.
For $1\leq j \leq r$, let $\ffdrj$ be the set of $f\in\ffdr$ such that 
$C_s$ has at most $j$ different values. We will be interested in
$j=1$ and 2. Clearly $\ffdro\subset\ffdrt\subset\cdots\subset \ffdr$
for each $\delta$ and $r$.

Let $h_C(x)\eed 
1+x'Cx$ for
$C\in \PP_d$ and $x\in\RR^d$. Then clearly $f\in\ffdro$ if and
only if for some 
$P\in\mm_{2r}$ and $C\in{\mathcal W}_{\delta}$
we have for all $x$
\begin{equation}\label{ffdroe}
f(x)\eee f_{P,C,r}(x)\eed
P(x)h_C(x)^{-r}.
\end{equation}

Each $C\in\PP_d$ can be diagonalized for an orthonormal basis
of eigenvectors with positive eigenvalues.
For $C=C_s\in{\mathcal W}_{\delta}$ as in the definition of $\ffdr$ we 
have for all $x\in\RR^d$ that $\delta|x|^2\leq x'Cx\leq 
|x|^2/\delta$. Set $\|f\|_{\sup}\eed\sup\{|f(x)|:\ x\in\RR^d\}$. 
Recall that $u\vee w\eed\max(u,w)$ 
and $u\wedge w\eed\min(u,w)$.


\begin{lem}\label{bounds}
For any 
$d=1,2,...,$ $r=1,2,...$,  $0<\delta<1$, 
and $f\in\ggdr$ we have
\fl
(a) $\|f\|_{\sup}\leq \delta^{-r}$,
\fl
(b) $\|f\|_{\sup}\geq (\delta/d)^{r}$.
\end{lem}
\pf
First let $f\in\ffdr$. For (a), we have for all $x$ that 
$$
|f(x)|\leq (1\vee |x|)^{2r}/(1+\delta|x|^2)^{r}
= (1\vee |x|^2)^{r}/(1+\delta|x|^2)^{r}\lee\delta^{-r},
$$
by considering the two cases $|x|\geq 1$ and $|x|<1$.
So $(a)$ follows in this case.

For $(b)$, for any $M$ with $1\leq M<\infty$ consider the point 
$x$ where $x_1=x_2=\cdots = x_d=M$. For this $x$ we get
$$
\|f\|_{\sup}\geq |f(x)| \geq M^{2r}/(1 + (dM^2/\delta))^{r}
= [\delta M^2/(\delta + dM^2)]^{r} \to (\delta/d)^r
$$
as $M\to +\infty$, proving $(b)$ for $f\in\ffdr$.

Now for any $f\in\ggdr$, we have $f\in{\mathcal F}_{\delta,v}$ for
some $v=1,...,r$, $\delta^{-v}\leq\delta^{-r}$ and
$(\delta/d)^v\geq (\delta/d)^r$, so the lemma follows.\qed
\vskip .1in
\begin{lem}\label{difflem}
For any 
$d=1,2,...$, $r=1,2,...$, and  $0<\delta<1$, 
let $f=f_{P,C,r}$ and $g=f_{P,D,r}$ for some 
$P\in\mm_{2r}$
and $C,D\in\PP_d$.
Then
\bge\label{difflemo}
(f-g)(x)\ \equiv\  {\frac{x'(D-C)xP(x)\sum^{r-1}_{j=0}h_D(x)^{r-1-j}
h_C(x)^j}{(h_Ch_D)(x)^{r}}}.
\end{equation}
For $1\leq k\leq l\leq d$ and $j=0,1,...,r-1$, let
\begin{eqnarray*}
h_{C,D,k,l,r,j}(x)&\eed & x_kx_lP(x)h_D(x)^{r-1-j}h_C(x)^j/(h_Ch_D)(x)^{r}\\
&= &x_kx_lP(x)h_C(x)^{j-r}h_D(x)^{-j-1}.
\end{eqnarray*}
Then each
$h_{C,D,k,l,r,j}$ is in ${\mathcal F}^{(2)}_{\delta,r+1}$ and
\bge\label{difflemt}
g-f\equiv -\sum_{1\leq k\leq l\leq d}
\sum_{j=0}^{r-1}(D_{kl}-C_{kl})(2-\delta_{kl})h_{C,D,k,l,r,j}.
\end{equation}
\end{lem}
\pf
Since $P(x)$ is a common factor we need only note that
$$
h_C(x)^{-r}-h_D(x)^{-r} = (h_D(x)^{r}-h_C(x)^{r})
/(h_Ch_D)(x)^{r}
$$
and expand the numerator by the identity
$$
U^{r}-V^{r} = (U-V)\sum^{r-1}_{j=0}U^{r-1-j}V^j
$$
for any two real numbers $U$ and $V$ to get (\ref{difflemo}).

The functions $h_{C,D,k,l,r,j}$ are clearly in 
${\mathcal F}^{(2)}_{\delta,r+1}$, and (\ref{difflemt}) follows
straightforwardly.
\qed
\vskip .1in
For 
$j=1,2,...,$ 
let $\ggdrj:= \bigcup^r_{v=1}\ffdvj$.
For any 
$f:\ \RR^d\to\RR$, define
$$
\|f\|^{\ast,j}_{\delta,r}
\eed \|f\|^{\ast,j}_{\delta,r,d}
$$
$$
\eed\inf\left\{\sumi_{s=1}|\lambda_s|: \ 
\exists g_s\in \ggdrj,\ 
s\geq 1,\ f\equiv \sumi_{s=1}\lambda_sg_s\right\},
$$
or $+\infty$ if no such $\lambda_s$, $g_s$ with 
$\sum_s|\lambda_s|<\infty$ exist. Lemma \ref{bounds}(a) implies
that for $\sum_s|\lambda_s|<\infty$ and $g_s\in\ggdr\equiv \ggdr^{(r)}$,
$\sum_s\lambda_sg_s$ converges absolutely and uniformly on 
$\RR^d$.
Let $Y^j_{\delta,r}
\eed Y^j_{\delta,r,d}$ be the set
of all functions $f$ from $\RR^d$ into $\RR$ such that
$\|f\|^{\ast,j}_{\delta,r}<\infty$. It's easily seen
that each $Y^j_{\delta,r}$ is a real vector space of functions
on $\RR^d$ and $\|\cdot\|^{\ast,j}_{\delta,r}$ is a seminorm on
it.
\begin{lem}\label{norms}
For any 
$j=1,2,...$, 
\fl
(a) If $f\in\ggdrj$ then $f\in Y^j_{\delta,r}$ and
$\|f\nrmdrstj \leq 1$.
\fl
(b) For any $g\in Y^j_{\delta,r}$, $\|g\|_{\sup}
\leq \|g\nrmdrstj/
\delta^{r}<\infty$.
\fl
(c) If $f\in\ggdrj$ then $\|f\nrmdrstj \geq (\delta^2/d)^r$.
\fl
(d) $\|\cdot\nrmdrstj$ 
is a norm on $Y^j_{\delta,r}$.
\fl
(e) $Y^j_{\delta,r}$ is complete for $\|\cdot\nrmdrstj$  
and thus a Banach space.
\end{lem}
\pf
Part $(a)$ is clear. For part $(b)$, let
$g=\sum_s\lambda_s g_s$ with $g_s\in
\ggdrj$
and apply Lemma \ref{bounds}$(a)$ to each $g_s$.
Part $(c)$ follows from part $(b)$ and
Lemma \ref{bounds}$(b)$.
For part $(d)$ we already noted that \anormdstj\ is a
seminorm. By part $(b)$, if $\|g\nrmdrstj=0$ then
$\|g\|_{\sup}=0$ so $g\equiv 0$, hence \anormdstj\ is
a norm.

For part $(e)$, let $\{f_k\}_{k\geq 1}$ be a Cauchy sequence in
$Y^j_{\delta,r}$ for \anormdstj. Then by part $(b)$ it is also
a Cauchy sequence for $\|\cdot\|_{\sup}$ and so converges
uniformly on $\RR^d$ to some function $f$. Taking a subsequence,
we get $f_{k_i}$ such that $\|f_{k_i}-f_k\nrmdrstj < 1/2^i$ for
all $k\geq k_i$ for $i=1,2,...$. Then the series
$f_{k_1}+\sumi_{i=1}f_{k_{i+1}}-f_{k_i}$ converges in
sup norm and in \anormdstj\ to $f$, writing
$f_{k_{i+1}}-f_{k_i}=\sum_s\lambda_{is}f_{is}$ for some
$f_{is}\in\ggdrj$ and $\sum_s|\lambda_{is}|\leq 1/2^i$, so
that $\sum_{i,s}|\lambda_{is}|\leq 1$.  It follows that 
$f\in Y^j_{\delta,r}$.
This finishes the proof.
\qed
\vskip .1in

\begin{lem}\label{inclusion}
For any 
$j=1,2,...$, we have 
$Y^j_{\delta,r}\subset Y^j_{\delta,r+1}$.
The inclusion linear map 
from $Y^j_{\delta,r}$ into $Y^j_{\delta,r+1}$ 
has norm at most 1.
\end{lem}
{\it Proof}. 
%
Simply 
$\ggdrj\subset {\mathcal G}^{(j)}_{\delta,r+1}$
and so $Y^j_{\delta,r}\subset Y^j_{\delta,r+1}$ with the inclusion
map having norm bounded above by 1.
\qed
\vskip .1in
\fl
{\bf Remark}. 
 Let spaces $Z^j_{\delta,r}$ be defined like spaces $Y^j_{\delta,r}$
but with $\ffdrj$ in place of $\ggdrj$.
Then the inclusion $Z^j_{\delta,r}\subset Z^j_{\delta,r+1}$  
doesn't hold, for example $x^2/(1+x^2)$ is in $Z^1_{\delta,1}$ 
but not in $Z^1_{\delta,2}$: suppose
$x^2/(1+x^2)=\sum_s \lambda_sx^4/(1+a_sx^2)^2$ with 
$\sum_s|\lambda_s|<\infty$ and $\delta<a_s<1/\delta$ for all $s$.
Dividing both sides by $x^2$ and letting $x\to 0$, the left side
approaches 1 and the right side 0.
\vskip .1in
In partial derivatives $\pard/\pard C_{kl}$ with respect to elements
of a symmetric matrix $C$, $C_{lk}\equiv C_{kl}$ will vary while
all other elements $C_{vw}$ are held fixed.

\begin{prop}\label{frechet}
Let $P\in\mmtr$, 
and consider the function $\phi(C,x)\eed  f_{P,C,r}(x)=
P(x)/h_C(x)^{r}$
from $\wwdelt\times \RR^d$ into $\RR$. Then:
\fl
(a) For each fixed $C\in\wwdelt$, $\phi(C,\cdot)\in \ffdro$.
\fl
(b) With respect to any entry $C_{kl}$ of $C$, $\phi(\cdot,x)$
has the partial derivative
$$
{\frac{\pard\phi(C,x)}{\pard C_{kl}}}\eee 
- {\frac{r(2-\delta_{kl})x_kx_lP(x)}{h_C(x)^{r+1}}}.
$$
\fl
(c) The map $C\mapsto \pard\phi(C,\cdot)/\pard C_{kl}$ is Lipschitz
from $\wwdelt$ into $Y^2_{\delta,r+2}$ for each $k,l=1,...,d$.
\fl
(d) The map $C\mapsto \phi(C,\cdot)$ from $\wwdelt$ into
$\ffdro\subset Y^1_{\delta,r}$, viewed as a map into the
larger space $Y^2_{\delta,r+2}$,
is Fr\'echet $C^1$.
\end{prop}
\pf
Part $(a)$ is clear, by (\ref{ffdroe}). Part $(b)$ follows by
elementary calculus. 

For part $(c)$, 
$Y^1_{\delta,r}\subset Y^2_{\delta,r}\subset Y^2_{\delta,r+2}$
where the first inclusion is immediate and the second follows
from Lemma \ref{inclusion} applied twice. From part (b),
the given map takes values in 
$Y^1_{\delta,r+1}\subset Y^1_{\delta,r+2}$ by Lemma \ref{inclusion}. 
It is Lipschitz into $Y^2_{\delta,r+2}$ by 
Lemma \ref{difflem} applied to $r+1$ in place of $r$.

For part $(d)$, and any $h\neq 0$ in $\RR$, let 
$C_{k,l,h}\in\wwdelt$ for $h$ small enough be the matrix which
equals $C$ except that $h$ is added to $C_{kl}$, and also to
$C_{lk}$ if $k\neq l$. Applying Lemma \ref{difflem},
we see that $[\phi(C_{k,l,h},\cdot)-\phi(C,\cdot)]/h
\in Y^2_{\delta,r+1}$ for $h$ small enough and equals,
for $D=C_{k,l,h}$, at any $x\in\RR^d$,
$$
-{\frac{(2-\delta_{kl})x_kx_lP(x)\sum^{r-1}_{j=0}h_D(x)^{r-1-j}
h_C(x)^j}{(h_Ch_D)(x)^{r}}}.
$$
It will be shown that as $D$ approaches $C$, i.e.\ $h\to 0$, the
function in the last display converges in $Y^2_{\delta,r+2}$ to the 
corresponding function with $D$ replaced by $C$. Considering one term 
at a time in the sum of $r$ terms, we get
$$
h_D(x)^{-j-1}h_C(x)^{j-r} - h_C(x)^{-r-1}
= h_C(x)^{j-r}[h_D(x)^{-j-1}-h_C(x)^{-j-1}].
$$
By the proof of Lemma \ref{difflem},
$$
h_D(x)^{-j-1}-h_C(x)^{-j-1}\ \equiv\  -h(2-\delta_{kl})x_kx_l
\sum^j_{i=0}h_C(x)^{i-j-1}h_D(x)^{-i-1}.
$$
Thus we get
$$
h(2-\delta_{kl})^2x_k^2x_l^2P(x)h_C(x)^{j-r}
\sum^j_{i=0}h_C(x)^{i-j-1}h_D(x)^{-i-1}\ \to\ 0
$$
in $Y^2_{\delta,r+2}$ as $h\to 0$. This implies
existence of the partial derivatives 
$$
\pard/\pard C_{kl}\left[C\mapsto \phi(C,\cdot)\in Y^2_{\delta,r+2}
\right].
$$
Part (c) gives their continuity. Continuous first partial derivatives
imply that the function $C\mapsto \phi(C,\cdot)$ is Fr\'echet
$C^1$ by known facts in analysis, completing the proof.
\qed
\vskip .1in

\begin{thm}\label{uniqueseries}
Let $r=1,2,...$, $d=1,2,...$, $0<\delta<1$, and $f\in Y^1_{\delta,r}$,
so that for some $a_s$ with $\sum_s |a_s|<\infty $ we have
$f(x)\equiv \sum_s a_sP_s(x)/(1+x'C_sx)^{k_s}$ for $x\in\RR^d$ where
each $P_s\in
{\MM\MM}_{2k_s}$, $k_s=1,...,r$, and $C_s\in\wwdelt$. Then $f$ can be
written as a sum of the same form in which the triples $(P_s,C_s,k_s)$
are all distinct. In that case, the $C_s$, $P_s$, $k_s$ and the 
coefficients $a_s$ are uniquely determined by $f$.
\end{thm}
\pf
Because $\sum |a_s|<\infty$, we can directly sum terms having
the same $P_s$, $C_s$ and $k_s$ into one. To prove the second conclusion
is equivalent to showing that then, if $f\equiv 0$ on $\RR^d$
we have all $a_s=0$.

Suppose the dimension $d=1$. Then each $C_s$ is a real number
with $\delta<C_s<1/\delta$, 
each $P_s(x)=x^{2k_s}$,
 and for each $s$ the denominator is $(1+C_sx^2)^{k_s}$. 
Each term in the sum is a rational function
of $x$ which we take as a rational function of a complex variable $z$.
Via a partial fraction decomposition (e.g.\ Knopp, 1947, Chapter 2),
we can rewrite the sum of at most $r$ terms
having a fixed value of $C_s=C$ as a sum
$$
\gamma_C+\sum_{k=1}^r {\frac{\alpha_{C,k}}{(z-z_C)^k}}
+ {\frac{\beta_{C,k}}{(z+z_C)^k}}
$$
where $z_C = i/\sqrt{C}$ and $\alpha_{C,k},\ \beta_{C,k}$ are
complex constants, $\gamma_C$ real.
It will suffice
to prove that 
\begin{equation}\label{sumabsfin}
\sum_s\left[|\gamma_{C_s}|+\sum_{k=1}^r|\alpha_{C_s,k}|+|\beta_{C_s,k}|
\right]\ <\ \infty,
\end{equation}
and given that, to show that for $C=C_s$ and $k=k_s$ for any $s$,
\begin{equation}\label{alphabetzer}
\alpha_{C,k}=\beta_{C,k}=0, 
\end{equation}
because if
we consider the largest value of $k=k_s$ with $a_s\neq 0$
and $C_s=C$, we will get $\alpha_{C,k}\neq 0$.

{\it Proof of (\ref{sumabsfin})}:
We will show that when $\delta<C<1/\delta$ and we take the
following partial fraction decomposition,
\begin{equation}\label{partfracC}
\frac{z^{2k}}{\left(1+Cz^2\right)^k}=\frac{1}{C^k}+
\sum_{v=1}^{k}\left(\frac{A_v}{\left(z-z_C\right)^v}
+\frac{B_v}{\left(z+z_C\right)^v}\right),
\end{equation}
then $\sum_v\left(\left|A_v\right|+\left|B_v\right|\right)$
has a suitable upper bound. In fact it will be shown in (\ref{avbvbd})
to be bounded above by $2^{k}\delta^{-3k/2}$,
and this will suffice.

Note that we have the term
$C^{-k}$ in the partial fraction decomposition (\ref{partfracC})
since this is the limit of the left-hand side as
$z\to\infty$. Let $a_{C,k}$ be the coefficient of $x^{2k}/h_C(x)^k$
in the statement of Theorem \ref{uniqueseries}.
Then 
$$
|\gamma_C|\lee \sum_{k=1}^rC^{-k}|a_{C,k}|\lee 
\delta^{-k}\sum_{k=1}^r|a_{C,k}|
$$
and $\sum_s|\gamma_{C_s}|<\infty$.

Using that $z_C=
{i}/{\sqrt{C}}$, we 
multiply both sides of (\ref{partfracC}) by $(-C)^k$,
then use $C = -z_C^{-2}$ to get
$$
\frac{
\left(
{z}/{z_C}\right)^{2k}}
{\left(1-\left(
{z}/{z_C}\right)^{2}\right)^k}=
(-1)^k
+\sum_{v=1}^{k}\left(\frac{A_v}{z_{C}^{2k+v}
\left(\frac{z}{z_C}-1\right)^v}+\frac{B_v}{z_C^{2k+v}
\left(\frac{z}{z_C}+1\right)^v}\right).
$$

Let now $A_v'=
{A_v}/{z_C^{2k+v}}$ and $B_v'=
{B_v}/{z_C^{2k+v}}$.
Note that $A_v'$ and $B_v'$ are the coefficients in the following
decomposition:
$$
\frac{\xi^{2k}}{\left(1-\xi^2\right)^k}
=(-1)^k+\sum_{v=1}^{k}\left(\frac{A_v'}{\left(\xi-1\right)^v}
+\frac{B_v'}{\left(\xi+1\right)^v}\right).
$$


To bound sums of absolute values of $A_v'$ and $B_v'$,
the following lemma will be proved:

\begin{lem}\label{specipartfrac}
(a) For any positive integers $m$ and $n$, and any complex $z\neq \pm 1$,
we have
\begin{equation}\label{oneinnum}
{\frac{1}{(1-z)^m(1+z)^n}} = 
\sum^m_{j=1}{\frac{a_j}{(1-z)^j}}+\sum^n_{i=1}{\frac{b_i}{(1+z)^i}}
\end{equation}
where $a_j\eed a_j^{m,n}>0$ and $b_i\eed b_i^{m,n}>0$ for all 
$j=1,...,m$ and $i=1,...,n$, and $\sum^m_{j=1}a_j+\sum^n_{i=1}
b_i=1$.
\fl
(b) We have for $k=0,1,...,m-1$
\begin{equation}\label{afmla}
a_{m-k}=\frac{1}{2^{n+k}}{n+k-1 \choose k}.
\end{equation}
\fl
(c) For $i=0,1,\dots,n-1$
\begin{equation}\label{bfmla}
b_{n-i}=\frac{1}{2^{m+i}}{m+i-1 \choose i}.
\end{equation}
\fl
(d) For each positive integer $k$ and $z\neq \pm 1$,
\begin{equation}\label{alphas}
\frac {z^{2k}}{(1-z^2)^k} = 
(-1)^k+\sum^k_{j=1}\alpha_j\left[{\frac{1}{(1-z)^j}}+{\frac{1}{(1+z)^j}}
\right]
\end{equation}
where $\alpha_j=\alpha_j^{(k)}$ are real numbers depending on $k$ 
and $1+2\sum^k_{j=1}|\alpha_j|\leq 2^k$.
\end{lem}
\pf
For (a), known facts about partial fraction decompositions provide 
a decomposition of the given form for some (real) coefficients
$a_j$ and $b_i$, which sum to $1$ by evaluation at $z=0$.
Positivity of $a_j$ and $b_i$ will follow from parts (b) and (c).

For (b), multiplying both sides of (\ref{oneinnum}) by $(z-1)^m$,
we get $(-1)^m(1+z)^{-n}$, which is holomorphic except at $z=-1$.
Comparing Taylor coefficients at $z=1$ of orders $0,1,...,m-1$ we
get (\ref{afmla}). For (c), a proof of (\ref{bfmla}) is symmetric.


Now for (d), clearly the left side has a partial fraction
decomposition with constant term $(-1)^k$, letting
$z\to+\infty$, and with coefficients times $(1-z)^{-j}$ 
and $(1+z)^{-j}$ for $j=1,...,k$.  The coefficients
of $(1-z)^{-j}$ and $(1+z)^{-j}$ are equal since
interchanging $z$ and $-z$ preserves the left side.
Taking a binomial expansion of
$z^{2k} = [(z^2-1)+1]^k$, the fraction on the left side equals
$
\sum^k_{j=0}{k\choose j}(-1)^j
(1-z^2)^{j-k}.
$
The sum of the absolute values of the coefficients in
the binomial expansion is $(1+1)^k = 2^k$. Applying part (a) to
each term $(1-z^2)^{j-k}$ for $j<k$ gives
a sum of the stated form. For $j=k$ we get
the constant term $(-1)^k$ as stated, and the bound
$2^k$ holds, proving (d) and the lemma.
\qed
\vskip .1in
It turns out that the upper bound $2^k$ in part 
(d) can be improved to $2\cdot(4/3)^k$.
Such a bound and its sharpness will follow from Theorem 
\ref{coefinequals} after the proof of Theorem \ref{uniqueseries}.

Returning to the proof of Theorem \ref{uniqueseries},
we have $A_v'=(-1)^v\alpha_v$ and $B_v'=(-1)^v\beta_v$ for each $v$,
so $|A_v'|=|\alpha_v|$ and $|B_v'|=|\beta_v|$.
Thus 
$$
\left|A_v\right|\le |\alpha_v|\left|z_C\right|^{2k+v}
\le |\alpha_v|\cdot\delta^{-3k/2}
$$
and, similarly, $\left|B_v\right|\le |\beta_v|\cdot\delta^{-3k/2}$. 
Finally,
we get that the coefficients in (\ref{partfracC}) are bounded by
\begin{equation}\label{avbvbd}
\delta^{-k}+\sum_v\left(\left|A_v\right|+\left|B_v\right|\right)
\le
{2^{k}}/{\delta^{3k/2}}.
\end{equation}
Clearly, ${2^{k}}/{\delta^{3k/2}}\leq {2^{r}}/{\delta^{3r/2}}$
for $k=1,...,r$, 
which finishes the proof of (\ref{sumabsfin}).

 So now, to prove (\ref{alphabetzer}), we can assume we have a sum 
$\sum_s a_s/(x-z_s)^{k_s}
= 0$ for all real $x$ where $k_s=1,...,r$, 
the pairs $(k_s,z_s)$ are different for different $s$, 
$\Re z_s=0$, $\delta\leq|\Im z_s|\leq 1/\delta$,
and $\sum_s|a_s|<\infty$.
The series converges uniformly
in any compact subset of the complement of $\{z:\ \Re z = 0,\ 
\delta\leq |\Im z|\leq 1/\delta\}$.

By analytic continuation, $g$ is holomorphic in the whole complement
and $0$ there. We will get a contradiction from that.



The following argument has been known apparently at least
since the 1920's. 
Ross and Shapiro (2002), Proposition 3.2.2 p.\ 18, give such an
argument for $r=1$, but it extends directly to any $r$.

Let $s$ be such that
$k_s$ is maximal, so we can assume $k_s=r$. It will be shown that
$t^rf(t+iy_s)\to a_s$ as real $t\to 0$, by dominated convergence for 
sums. We have
$|t^ra_v/(t + i(y_s-y_v))^{k_v})|\leq |a_v|$ for all $v$ since $k_v\leq 
r$, so we have
domination by a summable sequence. Moreover if $s\neq v$ the $v$ term 
approaches 0 as $t\to 0$
either because $y_s\neq y_v$ or they are equal and $k_v<r$. For $v=s$ we 
get $t^ra_s/t^r\equiv a_s\to a_s$ as $t\to 0$. By this contradiction,
the conclusion must hold for $d=1$.

%

Now 
let $d\geq 2$.
In this case there are still finitely
many possibilities for the polynomials $P_s$.
Call a set $A\subset\RR^d$ {\it algebraic} if it is of the
form $A=\{x:\ P(x)=0\}$ for some polynomial $P$, not
identically 0.  
It's easily seen that 
any countable union of algebraic sets in $\RR^d$ has dense
complement:
since an algebraic set $A$
is closed,
it suffices to show that $A$ is nowhere dense,
then apply the Baire category theorem. If $A$ were dense
in a non-empty open set $U$ it would include $U$.
So $P\equiv 0$ on $U$
and so $P$
is identically 0, a contradiction.

Now consider the cube $K_d$ of all vectors $\alpha = 
(\alpha_1,...,\alpha_d)$ with $1\leq\alpha_k\leq 2$
for all $k=1,...,d$. For each $\alpha\in K_d$, 
consider the function of one variable
$f_{\alpha}(t)\eed f(t\alpha)$, $t\in\RR$. For a dense
set of values of $\alpha\in K_d$, for all $s$, the numerator
of the $s$ term in $f_{\alpha}$ is a non-zero multiple of 
$t^{2k_s}$, and the numbers $\alpha'C_s\alpha$ are distinct
for different $C_s$, since these properties hold outside a
countable union of algebraic sets.
Choose and fix a value of $\alpha$ having
these properties. Then we can apply the $d=1$ case for
that value of $\alpha$, proving Theorem \ref{uniqueseries}.
\qed
\vskip .1in


To improve the bound in Lemma \ref{specipartfrac}(d), multiplying 
both sides of (\ref{alphas}) by $(1-z)^k$, we get
\begin{equation}\label{betteralph}
\left[\frac {z^{2}}{z+1}\right]^k = 
\sum^{k-1}_{v=0}\alpha_{k-v}(-1)^{v}(z-1)^v + g(z)(z-1)^k
\end{equation}
for some function $g$ holomorphic in a neighborhood of $z=1$.
It follows that $(-1)^{j}\alpha_{k-j}$ for $j=0,1,\dots,k-1$ are 
the Taylor coefficients of order $j$ of $[z^2/(z+1)]^k$ around
$z=1$ or equivalently, of $[(z+1)^2/(z+2)]^k =
\left (z + {\frac 1{z+2}}\right)^k$ around $z=0$. Since
the signs $(-1)^{j}$ don't affect the absolute values,
we need to bound the sum of the absolute values of these
Taylor coefficients through order $k-1$.
In Theorem \ref{coefinequals}, equation (\ref{thmpos}) we will see
that $(-1)^{j}\alpha_j^{(k)}>0$ for all $k=1,2,\dots$ and $j=0,1,\dots,
k-1$.

Let $c_i^n$ be the coefficients of the Taylor series of 
$\left(z+\frac{1}{2+z}\right)^n$ at $z=0$, namely
$$
\left(z+\frac{1}{z+2}\right)^n=c_{0}^n+c_{1}^nz+
\cdots+c_n^nz^n+c_{n+1}^nz^{n+1}+\cdots.
$$

Plugging in $z=1$, we get that
\begin{equation}\label{fthton}
\left(\frac{4}{3}\right)^n=c_0^n+c_1^n+\cdots+c_n^n+c_{n+1}^n+\cdots.
\end{equation}

In what follows, ${x\choose k}$ is defined (as usual) as
$x(x-1)\cdots(x-k+1)/k!$ for any real number $x$ and integer 
$k\geq 0$, and as $0$ otherwise (specifically, if $k<0$).

\begin{thm}\label{coefinequals}
For any integers $n\ge0$ and $l\ge0$
\beq\label{thmgen}
c_{n+l+1}^n=\frac{\left(-1\right)^l}{2^{2n+l+1}}\sum_{a=0}^l
{l \choose a}{2n \choose n-1-a}\left(-1\right)^a.
\end{equation}
We also have
\begin{equation}\label{thmuppersum}
T_n \eed\sum_{l=0}^{\infty}c_{n+1+l}^n\eee 
  \sum_{b=0}^{n-1}\frac{{2n \choose b}}{2^{2n}}\cdot\frac{3^b}{3^n}
> 0.
 \end{equation}

Also for any $0\le l\le n$
\begin{equation}\label{thmpos}
c_{n-l}^n=\frac{1}{2^{2n-l}}\sum_{a=l}^n{a \choose a-l}
{2n \choose n+a} > 0.
\end{equation}
In particular,
for any integer $n\ge 0$
\beq\label{thmspec}
c_n^n=\frac 12 + \frac 1{2^{2n+1}}{2n \choose n}.
\end{equation}
Let $S_n\eed \sum^{n-1}_{l=0}c^n_l=\sum^{n-1}_{l=0}|c^n_l|$. Then
\begin{equation}\label{Snbdasymp}
\left({\frac 43}\right)^n-{\frac 12} 
\geq S_n=
\left({\frac 43}\right)^n -{\frac 12}  - {\frac 1{\sqrt{\pi n}}}
+ o\left({\frac 1{\sqrt{n}}}\right).
\end{equation}
\end{thm}
\pf
In this proof, identities (\ref{vander1}) and (\ref{vander2})
can be traced back to a classical identity given by
Vandermonde (1772),
$\sum_{k\in\ZZ} {r \choose m+k} {s \choose n-k} = {r+s \choose m+n}$,
which is 
sometimes called Vandermonde's convolution. It holds for any
values of $r$ and $s$ provided that $m$ and $n$ are integers.
However, according to
Graham, Knuth and Patashnik (1994, \S5.1),  
this identity was known to Chu Shih-Chieh in China as
early as 1303.

 We will heavily use the following identity:
\beq\label{ident}
2^{-2n}{2n \choose n+a}\left(-1\right)^a=\sum_{b=0}^n{n \choose b}
{2b \choose b+a}2^{-2b}\left(-1\right)^{b},
\end{equation}
for which we don't know a reference,
so first let us show that it is valid. Note that we can rewrite it
as follows:
$$
2^{-2n}{2n \choose n+a}\left(-1\right)^{n+a}=\sum_{b=0}^n
{n \choose b}\left(-1\right)^{n-b}{2b \choose b+a}2^{-2b}.
$$
On the left is the coefficient of $z^{n+a}$ in 
the polynomial $\left(z-1\right)^{2n}2^{-2n}$. Let us compute that 
coefficient in a different way:
$$
\left(z-1\right)^{2n}2^{-2n}=\left(\frac{z^{2}-2z+1}{4}\right)^n=
\left(\left(-z\right)+\frac{(z+1)^{2}}{4}\right)^n.
$$
Using the binomial theorem twice, we have that:
\begin{eqnarray*}
\left(z-1\right)^{2n}2^{-2n} & = & \sum_{b=0}^n{n \choose b}
\left(-z\right)^{n-b}2^{-2b}(z+1)^{2b}\\
 & = & \sum_{b=0}^n{n \choose b}\left(-1\right)^{n-b}2^{-2b}z^{n-b}
\cdot\sum_{k=0}^{2b}{2b \choose k}z^k\\
 & = & \sum_{b=0}^n\sum_{k=0}^{2b}{n \choose b}{2b \choose k}
\left(-1\right)^{n-b}2^{-2b}z^{n-b+k}.
\end{eqnarray*}
Now in order to get the coefficient of $z^{n+a}$, we need to
plug in $k=b+a$, and this leads to 
(\ref{ident}).

Let us apply the binomial theorem to $\left(z+\frac{1}{z+2}\right)^n$.
We obtain that
\begin{eqnarray*}
\left(z+\frac{1}{z+2}\right)^n & = & \sum_{b=0}^n
{n \choose b}z^{n-b}\left(z+2\right)^{-b}\\
 & = & \sum_{b=0}^n{n \choose b}z^{n-b}\sum_{k\ge0}{-b \choose k}z^k2^{-b-k}\\
 & = & \sum_{b=0}^n\sum_{k\ge0}{n \choose b}{-b \choose k}2^{-b-k}z^{n-b+k}.
 \end{eqnarray*}
If $l=k-b-1$, then $n-b+k=n+l+1$, and 
we get a formula for the coefficient of $z^{n+l+1}$, namely
$$
c_{n+l+1}^n=\sum_{b=0}^n{n \choose b}{-b \choose b+l+1}2^{-2b-l-1},
$$
or, using the fact that 
${-b \choose b+l+1}={2b+l \choose b+l+1}\left(-1\right)^{b+l+1}$,
we get
\beq\label{cnnlo}
c_{n+l+1}^n=\sum_{b=0}^n{n \choose b}
{2b+l \choose b+l+1}\left(-1\right)^{b+l+1}2^{-2b-l-1}.
\end{equation}

This formula holds for any value of $l$, but it is not convenient
for our purposes. Let now $l$ be non-negative. Equating the 
coefficients of $z^{b+l+1}$ on both sides of
$
(z+1)^{2b+l}=(z+1)^l(z+1)^{2b},
$
we obtain that
\begin{equation}\label{vander1}
{2b+l \choose b+l+1}=\sum_{a=0}^l{l \choose a}{2b \choose b+a+1}.
\end{equation}

Now we rewrite 
(\ref{cnnlo})
as follows:

\begin{eqnarray*}
c_{n+l+1}^n
 & = & \sum_{b=0}^n{n \choose b}\left(-1\right)^{b+l+1}2^{-2b-l-1}
\sum_{a=0}^l{l \choose a}{2b \choose b+a+1}\\
 & = & (-1)^l\sum_{a=0}^l{l \choose a}2^{-l-1}\left[\,\sum_{b=0}^n
{n \choose b}{2b \choose b+a+1}2^{-2b}\left(-1\right)^{b+1}\right]\\
 & = & (-1)^l\sum_{a=0}^l{l \choose a}2^{-l-1}
\left[2^{-2n}{2n \choose n+a+1}\left(-1\right)^a\right]
 \end{eqnarray*}
using (\ref{ident}), which leads to 
(\ref{thmgen}).

Next, to prove (\ref{thmuppersum}),
we have from (\ref{thmgen}):
\begin{eqnarray*}
\sum_{l=0}^{\infty}c_{n+1+l}^n & = & 
\sum_{l=0}^{\infty}\frac{\left(-1\right)^l}{2^{2n+l+1}}
\sum_{a=0}^l{l \choose a}{2n \choose n-1-a}\left(-1\right)^a\\
 & = & \frac{1}{2^{2n+1}}\sum_{a\ge0}{2n \choose n-1-a}\left(-1\right)^a
\sum_{l\ge a}{l \choose a}\left(-2\right)^{-l}.
 \end{eqnarray*}

Note that if $x=-\frac{1}{2}$ then we have for any integer $a\geq 0$
\begin{eqnarray*}
\sum_{l\ge a}{l \choose a}\left(-2\right)^{-l} & = & 
\sum_{l\ge a}{l \choose a}x^l\\
 & = & \frac{x^a}{a!}
\sum_{l\ge a}l\left(l-1\right)\cdots\left(l-a+1\right)x^{l-a}\\
 & = & \frac{x^a}{a!}\sum_{l\ge a}\frac{d^a}{dx^a}x^l=
\frac{x^a}{a!}\frac{d^a}{dx^a}\sum_{l\ge0}x^l\\
 & = & \frac{x^a}{a!}\frac{d^a}{dx^a}\frac{1}{1-x}=\frac{x^a}{a!}
\frac{a!}{\left(1-x\right)^{a+1}}\\
 & = & \left(-1\right)^a\frac{2}{3^{a+1}}.
 \end{eqnarray*}
Here we added $1+x+\cdots+x^{a-1}$, but since its $a$th
derivative is zero, the equality still holds.
Thus,
\begin{eqnarray}\label{uppersum}
T_n= \sum_{l=0}^{\infty}c_{n+1+l}^n & = & \frac{1}{2^{2n}}
\sum_{a\ge0}{2n \choose n-1-a}3^{-a-1}\nonumber\\
 & = & \sum_{b=0}^{n-1}\frac{{2n \choose b}}{2^{2n}}\cdot\frac{3^b}{3^n},
 \end{eqnarray}
which is (\ref{thmuppersum}).


Next we want to prove (\ref{thmpos}).
Let $l>0$. Equating the coefficients of $z^{b-l+1}$ on both sides of
$
(z+1)^{2b-l}=(z+1)^{-l}(z+1)^{2b},
$
we obtain that
\begin{equation}\label{vander2}
{2b-l \choose b-l+1}=\sum_a{-l \choose a-l}{2b \choose b-a+1}.
\end{equation}
The sum in the last display is finite since ${-l \choose a-l}$ 
equals $0$ for $a<l$ and
${2b \choose b-a+1}$ is $0$ for $a>b+1$. Then equation 
(\ref{cnnlo}) can be transformed as follows:
\begin{eqnarray*}
c_{n-l+1}^n & = & \sum_{b=0}^n{n \choose b}
(-1)^{b-l+1}2^{-2b+l-1}\cdot
\sum_a{-l \choose a-l}{2b \choose b-a+1}\\
 & = & (-1)^l\sum_a{-l \choose a-l}2^{l-1}\left[
\,\sum_{b=0}^n{n \choose b}{2b \choose b+a-1}2^{-2b}
(-1)^{b+1}\right]\\
 & = & (-1)^l\sum_a{-l \choose a-l}2^{l-1}\left[{2n \choose n+a-1}
2^{-2n}(-1)^a\right]
\end{eqnarray*}
by (\ref{ident}).
Replacing ${-l \choose a-l}$ by ${a-1 \choose a-l}(-1)^{a-l}$, we
obtain:
$$
c_{n-l+1}^n=\sum_a{a-1 \choose a-l}{2n \choose n+a-1}2^{-2n+l-1},
$$
which leads to (\ref{thmpos}).

To prove (\ref{thmspec}), from (\ref{thmpos}) we have 
$
c_n^n=2^{-2n}\sum_{a=0}^n{2n \choose n+a},
$
so
\begin{eqnarray*}
2^{2n}c_n^n & = & \frac{1}{2}\sum_{a=0}^n\left({2n \choose n+a}+
 {2n \choose n-a}\right)\\
 & = & 
\frac{1}{2}\left({2n \choose n}+\sum_{a=0}^{2n}{2n \choose a}
\right)=\frac{2^{2n}+{2n \choose n}}{2},
\end{eqnarray*}
which leads to (\ref{thmspec}). 

Recall that
$a_n\sim b_n$ means $a_n/b_n\to 1$ as $n\to\infty$.
To study the behavior of $c_i^n$ for $n$ large, we have
the following,
which follows directly from Stirling's formula:
%
\beq\label{demoivre}
\frac{1}{2^{2n}}{2n \choose n}
\sim \frac 1{\sqrt{\pi n}}.
\end{equation}
This formula is a special case of the (local) central limit theorem for 
binomial probabilities first found by de Moivre (1733). 
Formula (\ref{demoivre}) immediately implies that
\begin{equation}\label{cnnfmla}
c_n^n-\frac{1}{2}\sim \frac{1}{2\sqrt{\pi n}}.
\end{equation}

Now in (\ref{thmuppersum}), $T_n= 
{\frac 1{3\cdot 4^n}}{{2n}\choose{n-1}}R_n$ where,
bounding by a geometric series,
$$
R_n\eed 1 + \frac{n-1}{3(n+2)} + {\frac {(n-1)(n-2)}{3^2(n+2)(n+3)}}
+ \cdots \lee {\frac 32}
$$
and $R_n\to 3/2$ as $n\to\infty$. It follows by Stirling's formula
that $T_n\sim 1/(2\sqrt{\pi n})$ as $n\goin$.

We have $c^n_l>0$ for $l=0,1,\dots,n-1$ by (\ref{thmpos}).
The bound on the left side of (\ref{Snbdasymp}) follows from
(\ref{fthton}), (\ref{thmuppersum}) and (\ref{thmspec}).
The asymptotic form on the right follows by using in addition
(\ref{cnnfmla}) and $T_n\sim 1/(2\sqrt{\pi n})$.
This finishes the proof of the theorem.
%
\qed
\vskip .1in

\vskip .1in
\fl
{\bf Remark}.
In the proof of Theorem \ref{uniqueseries} for $d=1$, one argument was
mentioned to be known according to Ross and Shapiro (2002, p.\ 18).
They say that the argument is closely related to an example of 
Poincar\'e (1883).
Ross and Shapiro don't give details of the argument (although it is 
short, as given in our proof), they just say the result 
``can be obtained from the dominated
convergence theorem.''  They also are treating points converging to 
the unit circle in
the plane along radii through 0, but the case of a line instead of the 
circle if anything seems easier.

Another interesting fact mentioned by Ross and Shapiro (2002) shows
that, again for $d=1$,
for such a proof, of uniqueness of $z_s$, $k_s$ and $a_s$ to
work, the condition that all the $z_s$ are on a line segment or
smooth curve is really used. In fact, there is a bounded sequence 
$\{z_j\}$ of distinct complex numbers with $|z_j|>1$ for all $j$,
such that the set of all limits of subsequences of $\{z_j\}$ is
exactly the unit circle $\{z:\ |z|=1\}$, and there exist
coefficients $a_j\in\CCC$ with $\sum_j|a_j|<\infty$, for
which the series $f(z)\eed \sum_j a_j/(z-z_j)$ converges
to $0$ for all $z$ with $|z|<1$, whereas on $|z|>1$, $f$ is
a meromorphic function with a pole of order 1 at each $z_j$.
Such an example is given in Ross and Shapiro (2002, Theorem
4.2.5), where in that theorem we take the special case
that $G$ is the open unit disk (just as in the proof given
by Ross and Shapiro) and $f\equiv 0$. Such examples
are said to have been first found by Wolff (1921).
\vskip .1in
As a corollary of Theorem \ref{uniqueseries},
in Lemma \ref{norms} in the special case $j=1$ we
get $\|f\|^{*,1}_{\delta,r} = 1$, showing that part (a) is
sharp and improving on part (c) in that case. 
\vskip .1in
\fl
{\bf Remark}. 
The equality is not true for general $j$, however.
Let $0<\delta<\alpha < 1$ and let 
$$
f(x) = x^4/[(1+\alpha x^2)(1+\alpha^{-1}x^2)]
$$
so that $f\in {\mathcal F}^{(2)}_{\delta,2}$. Then
$$
f(x)\equiv {\frac{\alpha}{1-\alpha^2}}\left(
{\frac{x^2}{1+\alpha x^2}} - {\frac{x^2}{1+\alpha^{-1}x^2}}
\right),
$$
so $\|f\|^{*,2}_{\delta,2}\leq 2\alpha/(1-\alpha^2)$.
This can be less than 1, in fact it's arbitrarily small
for small $\delta$ and $\alpha$.
\vskip .1in

In view of Theorem \ref{uniqueseries}, one might expect that if 
$$
\sum_s a_s/[(1+C_s x^2)^{u_s}(1+D_s x^2)^{v_s}]=0
$$
 with $\{C_s, D_s\} \neq \{C_t, D_t\}$ 
for $s\ne t$ then all $a_s = 0$. Unfortunately, this is wrong, as 
shown by the following simple example: 
$$
\frac 2{(1+x^2)(1+3x^2)}-\frac 
1{(1+x^2)(1+2x^2)} = \frac 1{(1+2x^2)(1+3x^2)}.
$$

For any $r=1,2,...$, $d=1,2,\dots$, 
any $P\in\MM\MM_{2r}$, and any $C\neq D$ in $\wwdelt$,
 let
$$
f_{P,C,D}(x)\eed f_{P,C,D,r}(x)\eed f_{P,C,D,r,d}(x)\eed 
$$
$$
\frac{P(x)}{(1+x'Cx)^{r}} - \frac{P(x)}{(1+x'Dx)^{r}}.
$$
By Lemma \ref{difflem},
for $C$ fixed 
and $D\to C$ we have $\|f_{P,C,D,r}\|^{*,2}_{\delta,r+1}\to 0$.
The following shows this is not true if $r+1$ in the norm is
replaced by $r$, even if the number of different $C_s$'s in
the denominator is allowed to be as large as possible, namely $r$:
\begin{prop}\label{multidneg}
For any $r=1,2,...$, $d=1,2,\dots$, $P\in\MM\MM_{2r}$,
and $C\neq D$ in $\wwdelt$, we have
$$
\|f_{P,C,D,r}\|^{*,r}_{\delta,r}=2.
$$
\end{prop}
\pf
First let $d=1$. Suppose that
\begin{equation}\label{fpcdexpr}
f_{P,C,D,r}(x)
\equiv f(x)\eed x^{2r}\sum_{s=1}^{\infty}\frac{a_s
}{(1
+x'C_{s1}x)\cdots(1+x'C_{sr}x)},
\end{equation}
where 
all $C_{sk}\in
(\delta,1/\delta)$,
and $\sum_s |a_s|<\infty$.
Then the series summing to $f$ in (\ref{fpcdexpr}) will converge
uniformly in all real $x$ and moreover for $x$ replaced
by complex $z$ with $|\Im z|
\leq\gamma$ for any
$\gamma$ with $0<\gamma<
\sqrt{\delta}$.   
Let $J\eed \{z\in\CCC:\ \Re z = 0,\ \delta\leq |\Im z|\leq 1/\delta\}.$
As seen previously, the series 
also will converge uniformly on any compact subset of the complement of
$J$.
A point $z_0 = iy$ in $J$ (so that $\delta\leq |y|\leq 1/\delta$
and $y$ is real) will be called a
{\it semipole of order} $r$ of a complex-valued function 
$g$ defined on $U\eed\CCC\setminus J$, where $r$ is a positive
integer, if for some $c\in\CCC$ with $c\neq 0$,
\begin{equation}\label{semipolelim}
\lim_{t\downarrow 0}t^rg(iy+t)\ =\ c.
\end{equation}

Let
$$
F(z)\ :=\ z^{2r}\left[\frac 1{(1+Cz^2)^r} - \frac 1{(1+Dz^2)^r}
-\sum^{\infty}_{s=1}\frac {a_s}{\prod^r_{k=1}(1+C_{sk}z^2)}\right],
$$
convergent to 0 on the complement of $J$ and uniformly on compact subsets
of the complement. Let $z:=x + i/\sqrt{C}$ and $x\to 0$, and consider
the behavior of $x^{r}F(z)$. For any $B>\delta$ such as $B=C_{sk}$, $C$ or
$D$, we will have $1+Bz^2= 1-B(C^{-1}-x^2) + 2iBx/\sqrt{C}$, which has
absolute value at least $2\delta |x|/\sqrt{C}$. Thus
$|x/(1+Bz^2)|\leq \delta^{-3/2}$ and for the limit we have dominated
convergence for sums. Moreover, the limit is 0 except for terms with
$C_{sk}=C$ for all $k=1,...,r$. Combining the $a_s$ for such terms into
one, we must have $a_s=1$ for such an $s$. We can do likewise for $D$ 
in place of $C$, which finishes the proof for $d=1$.

For $d\geq 2$, let $S^{d-1}$ denote the unit
sphere centered at the origin in $\RR^d$. We will
use the Haar measure (rotationally invariant Borel probability
measure) $\mu$ on $S^{d-1}$.
For any $e\in S^{d-1}$, define a function $g_e$ of
a real variable $t$ by
$$
g_e(t)\eed f(te)
\eee t^{2r}\sum_{s=1}^{\infty}
\frac{a_sP_s(e)}{\left(1+t^2(e'C_{s1}e)\right)
\cdots\left(1+t^2(e'C_{sr}e)\right)}.
$$

For any fixed matrix $M$, e.g.\ $M=C_{sk}$, let $A_{M,C}$
be the set of all $e$ in $S^{d-1}$ for which 
$e'Ce=e'Me$, or
$
e'(C-M)e=0.
$
If $C\neq M$, then not all eigenvalues of $C-M$ are zero, which
entails that the codimension of $A_{M,C}$ is at least $1$. Thus, we 
obtain that $\mu(A_{M,C})=0$ provided that $M\neq C$,
and likewise $\mu(A_{M,D})=0$ for $M\neq D$.

For each $e$ such that $P(e)\neq 0$ and $e'(C-D)e\neq 0$, as are true 
for almost all $e$,
the function $f_{e,P,C,D}(t)\eed f_{P,C,D}(te)$ has two poles of
order $r$ at the points $z_{e}^\pm=\frac{\pm i}{\sqrt{e'Ce}}$.
Then, $g_e$ must have semipoles of order $r$ at the same two points,
as shown in the $d=1$ case.
The same is true for the points
$\zeta_{e}^\pm=\frac{\pm i}{\sqrt{e'De}}$. 

For $\mu$-almost all $e$, we have $e'C_{sj}e\neq e'Ce$
for all $s$ and $j$ such that $C_{sj}\neq C$, and $P(e)\neq 0$.
For any such $e$, contributions to the limit $c$ at
$iy=z_e^{\pm}$, via dominated convergence for sums as in
the proof of Theorem \ref{uniqueseries}, can come only from
terms in (\ref{fpcdexpr}) in which $C_{sj}=C$ for all
$j=1,...,r$.  
Let $MM(r)\eed MM(r,d)$ be the finite
number of different monomials in ${\MM\MM}_{r,d}$, 
namely
$MM(r,d) = {{2r+d-1}\choose {d-1}}$, and call them
$Q_1,Q_2,\dots,Q_{MM(r)}$. Let $b_j=a_s$ when
$P_s=Q_j$ and $C_{si}=C$ for all $i=1,...,r$. Some of these
coefficients may be $0$. We then get an equation
$$
P(e) = \sum^{MM(r)}_{j=1}b_jQ_j(e)
$$
holding for $\mu$-almost all $e$, by applying (\ref{semipolelim})
to both sides of (\ref{fpcdexpr}) and using dominated convergence
for sums, as in the proof of Theorem \ref{uniqueseries}, or the
proof for $d=1$.
Multiplying both sides by $t^{2r}$ we get a polynomial equality
holding almost everywhere on $\RR^d$ and so everywhere. Since
different monomials are linearly independent, it follows that
$b_j=1$ for the $j$ such that $P=Q_j$ and other $b_j=0$.
The same argument with $D$ in place of $C$ gives that
 $a_s=-1$ for another coefficient in (\ref{fpcdexpr}) and so shows
that
$
\|f_{P,C,D,r}\|^{*,r}_{\delta,r}\geq 2.
$
From the original expression of $f_{P,C,D,r}$ it's clear that
the norm is $\leq 2$, so it follows that it equals $2$.


Here, in (\ref{fpcdexpr}) we could also consider on the right
terms where the numerators and denominators were of the same
form but of total degree $2r-2,2r-4,...,0$ rather than $2r$.
There would be no essential change in the proof since such
terms do not affect or contribute semipoles of order $r$.
The proposition is proved.
\qed
\vskip .1in


%
\vskip .1in
\fl
{\bf Remark}. 
In Lemma \ref{bounds} we had two-sided bounds
in the sup norm for functions in $\ggdr$ and in Lemma
\ref{norms}(b), the sup norm is bounded above in terms of
the $\|\cdot\|^{*,j}_{\delta,r}$ norm. There is no such
bound in the reverse direction and in fact, for example,
$Y^1_{\delta,r}$ is not complete in the sup norm in dimension 
$d=1$. To see this, for a fixed $\delta\in (0,1)$ let $
a\eed\delta$ and $b\eed 1/\delta$. For $n=1,2,...,$ let 
$C_{jn} \eed a + (j(b-a)/n) - 1/(2n)$ for $j=1,...,n$.
Then $a<C_{jn}<b$ for each $j$. Let
$f_n(x) \eed \sum^n_{j=1} x^{2r}/[n\{1+(C_{jn}x^2)^r\}]$.
Then it is not hard to show that $f_n$ converges in the sup
norm as $n\goin$ to the 
function $g(x)\eed x^{2r}\int_a^b dC/[(b-a)(1+Cx^2)^r]$,
but $g$ is not in $Y^1_{\delta,r}$, because $g$ has no
semipoles.

\vskip .1in
\centerline{REFERENCES}
\vskip .1in
\bbe
Archibald, R.\ C.\ (1926), A rare pamphlet of Moivre and some of his
discoveries, {\it Isis} {\bf 8}, 671-676.
\bbe
Graham, R.\ L., Knuth, D.\ E., and Patashnik, O.\ (1989, 2d ed.\ 1994),  
{\it Concrete Mathematics: A Foundation for Computer Science},
Ad\-di\-son-Wesley, Reading, MA.
\bbe
Knopp, K.\  {\it Theory of Functions}, Part 2, transl.\ from the
4th German ed.\ by F.\ Bagemihl, Dover, New York, 1947.
\bbe
de Moivre, A.\ (1733), Approximatio ad Summam Terminorum Binomii
$(a+b)^n$ in Seriem Expansii, 7 pp., London. Transl.\ by the
author and included in de Moivre (1756). Latin orig.\ repr.\
in Archibald (1926).
\bbe
de Moivre, A.\ (1756, posth.), {\it The Doctrine of Chances},
3d ed., London. Repr.\ F.\ Case, London, and Chelsea, New York,
1967.
\bbe
Poincar\'e, H.\ (1883), Sur les fonctions \`a espaces
lacunaires, {\it Acta Soc.\ Scient.\ Fennicae} {\bf 12},
341-350. [Reference according to Ross and Shapiro.]
\bbe
Ross, W.\ T., and Shapiro, H.\ S.\ (2002), 
{\it Generalized Analytic Continuation}, {\it University Lecture
Series} {\bf 25}, Amer.\ Math.\ Soc.
\bbe
Vandermonde, 
A.\ (1772),
``M\'{e}moire sur des irrationnelles de
diff\'{e}rens ordres avec une application au cercle,'' {\it Histoire de
l'Acad.\ 
Royale des Sciences} 
(1772), 
part 1, 71-72; {\it M\'{e}moires
Math.\
Physique, 
Tir\'{e}s des Registres de
l'Acad.\
Royale des Sciences} (1772), 489-498. [Reference according
to Graham et al.]
\bbe
Wolff, J.\ (1921), Sur les s\'eries $\sum A_k/(z-z_k)$,
{\it Comptes Rendus Acad.\ Sci.\ Paris} {\bf 173},
1057-1058,1327-1328.

\end{document}